\documentclass[number,citesort,secthm,seceqn,dvips]{arxbj}
\usepackage{upgreek}

\aid{0}
\volume{17}
\issue{1}
\pubyear{2011}
\firstpage{114}
\lastpage{137}
\doi{10.3150/10-BEJ263}

\makeatletter
\newcommand{\eqref}[1]{(\ref{#1})}
\newcommand{\ds}{\displaystyle}
\renewcommand{\citet}[1]{\cite{#1}}

\newremark{rem}{Remark}[section]

\makeatother

\begin{document}
\begin{frontmatter}

\title{On a fractional linear birth--death process}
\runtitle{Fractional birth--death process}

\begin{aug}
\author[a]{\fnms{Enzo} \snm{Orsingher}\corref{}\thanksref{e1}\ead[label=e1,mark]{enzo.orsingher@uniroma1.it}} \and
\author[a]{\fnms{Federico} \snm{Polito}\thanksref{e2}\ead[label=e2,mark]{federico.polito@uniroma1.it}}
\runauthor{E. Orsingher and F. Polito}
\address[a]{Department of Statistics, Probability and Applied Statistics,
Sapienza University of Rome, pl. A. Moro 5, 00185 Rome, Italy. \printead{e1,e2}}
\end{aug}

\received{\smonth{9} \syear{2009}}
\revised{\smonth{1} \syear{2010}}

%
\begin{abstract}
In this paper, we introduce and examine a fractional linear
birth--death process
$N_\nu ( t )$, $t>0$, whose fractionality is obtained by replacing
the time derivative with a fractional derivative in the system of
difference-differential
equations governing the state probabilities $p_k^\nu ( t )$, $t>0$,
$k \geq0$. We present a subordination relationship connecting $N_\nu
( t )$,
$t>0$, with the classical birth--death process $N ( t )$, $t>0$, by
means of
the time process $T_{2 \nu} ( t )$, $t>0$, whose distribution is related
to a time-fractional diffusion equation.

We obtain explicit formulas for the extinction probability $p_0^\nu (
t )$
and the state probabilities $p_k^\nu ( t )$, $t>0$, $k \geq1$, in the
three relevant cases $\lambda> \mu$, $\lambda< \mu$, $\lambda= \mu$
(where $\lambda$ and $\mu$ are, respectively, the birth and death
rates) and
discuss their behaviour in specific situations. We highlight the
connection of the fractional
linear birth--death process with the fractional pure birth process.
Finally, the mean values $\mathbb{E} N_\nu ( t )$ and
$\operatorname{\mathbb{V}ar}N_\nu ( t )$ are derived and analyzed.
\end{abstract}

%

\end{frontmatter}

\section{Introduction}
In a previous paper \cite{pol}, we constructed a fractional version of
the pure birth process
$\mathcal{N}_\nu ( t ), t>0$ (both in the general and in the linear case
denoted here as $M_\nu ( t ), t>0$),
by considering the fractional equations governing their distributions.
In this work, we examine the linear birth--death process $N_\nu ( t )$,
$t>0$, where the state probabilities
%
\begin{equation}
p_k^\nu ( t ) = \Pr \{ N_\nu ( t ) = k | N_\nu ( 0 ) = 1 \}
\end{equation}
are assumed to satisfy the fractional difference-differential equations
%
\begin{eqnarray}
\label{state-eq}
\frac{\mathrm{d}^\nu p_k ( t )}{\mathrm{d} t^\nu} &=& - ( \lambda+ \mu
) k p_k
( t ) + \lambda ( k-1 ) p_{k-1} ( t )\nonumber\\[-8pt]\\[-8pt]
&&{}+ \mu ( k+1 ) p_{k+1} ( t ),\qquad k \geq1,
0<\nu\leq1 .\nonumber
\end{eqnarray}
The fractional operator appearing in \eqref{state-eq} is defined as
%
\begin{equation}
\label{eq:fractionalDerivative}
\cases{
\ds \dfrac{\mathrm{d}^\nu f ( t ) }{\mathrm{d}t^\nu} = \dfrac{1}{\Gamma ( 1- \nu )}
\int_0^t \frac{ (\mathrm{d}/{\mathrm{d}s}) f
( s )}{ ( t-s )^\nu} \,\mathrm{d}s,
&\quad $0 < \nu< 1$, \cr
f' ( t ), &\quad $\nu=1$.
}
\end{equation}
The derivative \eqref{eq:fractionalDerivative} is usually called a
Caputo or
Dzherbashyan--Caputo fractional derivative and differs from the
classical Riemann--Liouville
derivative by exchanging the integral and derivative operators (see
\cite{podlubny}).
An advantage of Caputo over Riemann--Liouville is that
Caputo does not require fractional-order derivatives in the initial
conditions, which is good for practical purposes.
The positive parameters $\lambda$ and $\mu$ are, respectively, the
birth and death rates.

The exact distribution of the linear birth--death process reads (see
\citet{bailey}, page
91, \citet{feller1}, page~454)
%
\begin{equation}
\label{class-state}
p_k^1 ( t ) = ( \lambda- \mu )^2 \mathrm{e}^{- ( \lambda- \mu
)t} \frac{ ( \lambda ( 1-\mathrm{e}^{- ( \lambda- \mu ) t}
) )^{k-1}}{ ( \lambda- \mu \mathrm{e}^{- ( \lambda- \mu ) t}
)^{k+1}},\qquad  k\geq1, t>0, \mu\neq\lambda.
\end{equation}
When $\lambda=\mu$, the distribution \eqref{class-state} is much
simpler and takes the
form
%
\begin{equation}
\label{class-state-2}
p_k^1 ( t ) = \frac{ ( \lambda t )^{k-1}}{ (
1+ \lambda t )^{k+1}},\qquad  t>0, k \geq1 .
\end{equation}
The exact expressions for the extinction probabilities are
%
\begin{equation}
\label{class-ext}
p_0^1 ( t ) = \cases{
\dfrac{\lambda t}{1+ \lambda t}, & \quad $\lambda= \mu$, \cr
\dfrac{\mu- \mu \mathrm{e}^{-t ( \lambda- \mu ) }}{
\lambda- \mu \mathrm{e}^{- t ( \lambda- \mu ) }}, &\quad  $\lambda\neq\mu$.
}
\end{equation}
From \eqref{state-eq}, we can easily infer that the probability
generating function
of $N_\nu ( t ), t> 0$,
%
\begin{equation}
G_\nu ( u,t ) = \mathbb{E} u^{N_\nu ( t )},
\qquad
|u| \leq1, 0< \nu\leq1, t>0,
\end{equation}
satisfies the Cauchy problem
%
\begin{equation}
\label{eq:fractDiffEqGF}
\cases{
\dfrac{\partial^\nu}{\partial t^\nu} G_\nu ( u,t ) =
( \lambda u - \mu )
( u -1 ) \dfrac{\partial}{\partial u} G_\nu ( u,t ),
&\quad  $\nu\in(0,1], |u| \leq1$, \vspace*{2pt}\cr
G_\nu ( u,0 ) = u.
}
\end{equation}
We will show below that from \eqref{eq:fractDiffEqGF}, one can arrive
at the subordination
relationship
%
\begin{equation}
\label{sub-repr}
N_\nu ( t ) \stackrel{\mathrm{i.d.}}{=} N ( T_{2 \nu} ( t )
),\qquad  t>0,
\end{equation}
where $T_{2 \nu} ( t )$, $t>0$, is the random time process whose
distribution is
obtained by folding the solution of the following fractional diffusion equation:\vadjust{\eject}
%
\begin{equation}
\label{fractDiffusion}
\cases{
\dfrac{\partial^{2 \nu} q}{\partial t^{2 \nu}} =
\dfrac{\partial^2 q}{\partial x^2}, &\quad
$0 < \nu\leq1$, $x \in\mathbb{R}$, $t>0$, \vspace*{2pt}\cr
q ( x,0 ) = \delta ( x ).
}
\end{equation}
The process $N ( t ), t>0$, found in \eqref{sub-repr}, is the classical
linear birth--death process
whose distribution is given in \eqref{class-state}, \eqref
{class-state-2} and \eqref{class-ext}.
A relationship similar to \eqref{sub-repr} also holds for the
fractional pure birth process
\citet{pol} and the fractional Poisson process \citet{fr-poisson}.
In this context, it represents the main tool of our analysis and leads
to a number of
interesting explicit distributions.\vadjust{\goodbreak} We consider the subordinator
related to
\eqref{fractDiffusion} because the probability generating function of
the distribution of \eqref{sub-repr}
satisfies the simplest fractional equation generalizing the classical one.

For the extinction probabilities of the fractional linear birth--death
process, we have the following
attractive formulas:
%
\begin{equation}
\label{total-extinction}
p_0^\nu ( t ) = \cases{
\ds\frac{\mu}{\lambda}
- \frac{\lambda- \mu}{\lambda}
\sum_{m=1}^{+ \infty} \biggl( \frac{\mu}{\lambda} \biggr)^m E_{\nu,1}
\bigl( -t^\nu ( \lambda- \mu
) m \bigr), &\quad $ \lambda>\mu$, \cr
\ds1 - \frac{\mu- \lambda}{\lambda}
\sum_{m=1}^{+\infty} \biggl( \frac{\lambda}{\mu} \biggr)^m E_{\nu,1}
\bigl( - t^\nu ( \mu- \lambda ) m \bigr), & \quad $\lambda<\mu$, \cr
\ds1 - \int_0^{+\infty} \mathrm{e}^{-w} E_{\nu,1} ( - \lambda t^\nu w ) \,\mathrm{d}w,
&\quad $ \lambda=\mu$
}
\end{equation}
for $t>0, 0< \nu\leq1$.

The function $E_{\alpha,\beta} ( x )$ appearing in \eqref{total-extinction}
is the generalized Mittag--Leffler function, defined as
%
\begin{equation}
\label{mittag}
E_{\alpha,\beta} ( x ) = \sum_{m=0}^{+ \infty} \frac{x^m}{
\Gamma ( \alpha m + \beta ) },\qquad  x \in\mathbb{R}, \alpha
> 0, \beta> 0 .
\end{equation}

From \eqref{total-extinction}, we can easily retrieve the classical
extinction probabilities
\eqref{class-ext} for $\nu=1$ by keeping in mind that $E_{1,1} ( x ) = \mathrm{e}^x$.

For the state distributions $p_k^\nu ( t )$, $t>0$, $k \geq1$,
we have formulas similar to
\eqref{total-extinction}, but with a more complicated structure:
%
\begin{eqnarray}
\label{state}
\hspace*{-22pt}
p_k^\nu ( t ) &=& \cases{
\ds\biggl( \frac{\lambda-\mu}{\lambda} \biggr)^2
\sum_{l=0}^\infty\pmatrix{l+k\cr l} \biggl( \frac{\mu}{\lambda} \biggr)^l
\sum_{r=0}^{k-1}
( -1 )^r \pmatrix{k-1\cr r}\cr
\ds\hphantom{\biggl( \frac{\lambda-\mu}{\lambda} \biggr)^2
\sum_{l=0}^\infty\pmatrix{l+k\cr l} \biggl( \frac{\mu}{\lambda} \biggr)^l
\sum_{r=0}^{k-1}} {}\times E_{\nu,1} \bigl( - ( l+r+1 ) ( \lambda-\mu ) t^\nu
\bigr), &\quad $ \lambda>\mu$, \cr
\ds\biggl( \frac{\lambda}{\mu} \biggr)^{k-1}
\biggl( \frac{\mu-\lambda}{\mu} \biggr)^2
\sum_{l=0}^\infty
\pmatrix{l+k\cr l} \biggl( \frac{\lambda}{\mu} \biggr)^l\cr
\ds\hphantom{\biggl( \frac{\lambda}{\mu} \biggr)^{k-1}
\biggl( \frac{\mu-\lambda}{\mu} \biggr)^2
\sum_{l=0}^\infty}{}\times\sum_{r=0}^{k-1} ( -1 )^r
\pmatrix{k-1\cr r}\cr
\ds\hphantom{\biggl( \frac{\lambda}{\mu} \biggr)^{k-1}
\biggl( \frac{\mu-\lambda}{\mu} \biggr)^2
\sum_{l=0}^\infty\times\sum_{r=0}^{k-1} }{}\times
E_{\nu,1} \bigl( - ( l+r+1 ) (\mu- \lambda
) t^\nu \bigr), &\quad $ \lambda<\mu$, \cr
\ds\frac{
( -1 )^{k-1}
\lambda^{k-1}}{k!} \frac{\mathrm{d}^k}{\mathrm{d} \lambda^k} \bigl[ \lambda \bigl( 1-
p_0^\nu ( t ) \bigr) \bigr], &\quad $ \lambda=\mu$.
}
\end{eqnarray}
Also from \eqref{state}, for $\nu=1$, one can reobtain the
distributions \eqref{class-state} and
\eqref{class-state-2}.

We will show below that the probabilities $p_k^\nu ( t )$, $t>0$, $k
\geq1$, appearing
in \eqref{state} are strictly related to the distributions of the
fractional linear pure
birth process $M_\nu(t)$, $t>0$, with an arbitrary number of
progenitors and a birth rate equal to $\lambda-\mu$
with $\lambda>\mu$. In particular, we can extract from the first line
of \eqref{state}
that
%
\begin{eqnarray}
&& \Pr \{ N_\nu ( t ) = k | N_\nu ( 0
) = 1 \}\nonumber \\
&&\quad  = \frac{\lambda- \mu}{\lambda} \sum_{l=0}^\infty \biggl[
\biggl( 1+\frac{\mu}{k ( \lambda-\mu )} \biggr)
\Pr\{ \mathcal{G}=l \} + \frac{\mu}{k} \frac{\mathrm{d}}{\mathrm{d} \mu} \Pr
\{ \mathcal{G}=l \}\biggr] \\
&&\qquad {}\times\Pr \{ M_\nu ( t ) = k+l |
M_\nu ( 0 ) = l+1 \}, \nonumber
\end{eqnarray}
where
%
\begin{equation}
\Pr \{ \mathcal{G}=l \} = \biggl( 1 -\frac{\mu}{\lambda} \biggr)
\biggl( \frac{\mu}{\lambda} \biggr)^l,\qquad  l \geq0,
\end{equation}
is a geometric law for the number of progenitors.
We also note that for $\lambda=\mu$, the distribution~\eqref{state} can be expressed in terms of the extinction
probability \eqref{total-extinction} by means of
%
\begin{equation}
\Pr \{ N_\nu ( t ) = k \} =
\frac{ ( -1 )^{k-1} \lambda^{k-1}}{k!} \frac{\mathrm{d}^k}{\mathrm{d} \lambda^k}
\bigl[ \lambda \bigl( 1 - \Pr \{ N_\nu ( t ) = 0 \}
\bigr) \bigr],\qquad  k \geq1, t>0 .
\end{equation}
The extinction probability \eqref{total-extinction} can be viewed as
being a suitable
weighted mean of the waiting times of the fractional Poisson process
$\mathcal{P}_\lambda^\nu ( t
), t>0$, for which it is well known that \citet{fr-poisson}
%
\begin{equation}
\Pr \{ \mathcal{N}_{\mathcal{P}_\lambda^\nu} ( t ) = 0 \}
= E_{\nu,1} ( - \lambda t^\nu ),\qquad  t>0 .
\end{equation}

The fractional linear birth--death process dealt with in this paper
provides a generalization
of the classical linear birth--death process and may well prove to be
capable of modeling queues in
service systems, epidemics and the evolution of populations under
accelerating conditions.
The introduction of the fractional derivative furnishes the system
with a global memory.
Furthermore, the qualitative features illustrated in the last section show
that the fractional counterpart of the linear birth--death process has a
faster mean evolution (and variance expansion), as was pointed
out in similar fractional generalizations, for example, for the Poisson
process (see
\citet{cahoy,sibatov,laskin,fr-poisson}), for
fractional branching processes \citet{sib} and for pure birth processes
\citet{pol}.

\section{The extinction probability of the fractional linear
birth--death process}
\label{sec}
We begin this section by proving the subordination relationship \eqref{sub-repr}
which is relevant to all of the distributional results of this paper.
\begin{thm}
\label{teo-ite}
The fractional linear birth--death process $N_\nu ( t ), t>0$,
can be represented as
%
\begin{equation}
\label{ite}
N_\nu ( t ) = N ( T_{2 \nu} ( t ) ),\qquad
t>0, 0<\nu\leq1,
\end{equation}
where $N ( t )$, $t>0$, is the classical linear birth--death process
and $T_{2 \nu} ( t ), t>0$, is a random process whose
one-dimensional distribution coincides with the folded solution
of the fractional diffusion equation
%
\begin{equation}
\label{funda}
\dfrac{\partial^{2 \nu} q}{\partial t^{2 \nu}} =
\dfrac{\partial^2 q}{\partial x^2},\qquad  0 < \nu\leq1, x \in\mathbb{R},
t>0,
\end{equation}
subject to the initial conditions $q ( x,0 ) = \delta ( x )$
for $0<\nu\leq1$ and also $q_t ( x,0 ) = 0$ for $1/2 < \break\nu\leq1$.
\end{thm}

\begin{pf}
The Laplace transform $\tilde{G}_\nu ( u,z )
= \int_0^\infty \mathrm{e}^{-zt} G_\nu ( u,t ) \,\mathrm{d}t$,
applied to the fractional PDE
%
\begin{equation}
\label{genf}
\cases{
\dfrac{\partial^\nu}{\partial t^\nu} G_\nu ( u,t ) =
( \lambda u - \mu )
( u -1 ) \dfrac{\partial}{\partial u}
G_\nu ( u,t ),
&\quad $ 0 < \nu\leq1$,\vspace*{2pt} \cr
G_\nu ( u,0 ) = u,
}
\end{equation}
yields
%
\begin{equation}
\label{gtilde}
z^\nu\tilde{G}_\nu ( u,z )
- z^{\nu-1} u = ( \lambda u-\mu )
( u-1 ) \frac{\partial}{
\partial u} \tilde{G}_\nu ( u,z ),\qquad
0<\nu\leq1, z>0, |u| \leq1 .
\end{equation}
We now observe that
%
\begin{equation}
\tilde{G}_\nu ( u,z ) =
\int_0^\infty \mathrm{e}^{-zt} \Biggl[ \sum_{k=0}^\infty u^k
\Pr \{ N_\nu ( t ) = k \} \Biggr] \,\mathrm{d}t.
\end{equation}
If \eqref{ite} holds, then
%
\begin{eqnarray}
\label{iterg}
\tilde{G}_\nu ( u,z ) & =&
\int_0^\infty \mathrm{e}^{-zt} \Biggl[ \sum_{k=0}^\infty u^k
\int_0^\infty\Pr \{ N ( s ) = k \}
\Pr \{ T_{2 \nu} ( t ) \in \mathrm{d}s \} \Biggr] \,\mathrm{d}t\nonumber \\
& = &\int_0^\infty \mathrm{e}^{-zt} \Biggl[ \int_0^\infty G ( u,s )
\Pr \{ T_{2 \nu} ( t ) \in \mathrm{d}s \} \Biggr] \,\mathrm{d}t
\\
& =& \int_0^\infty G ( u,s ) z^{\nu-1} \mathrm{e}^{-sz^\nu} \,\mathrm{d}s .
\nonumber
\end{eqnarray}
In the last step, we applied the folded version of
equation (3.3) in \citet{ors2004} for $c=1$, that being therefore
%
\begin{equation}
\int_0^\infty \mathrm{e}^{-zt} \Pr \{ T_{2 \nu} ( t )
\in \mathrm{d}s \} = \mathrm{e}^{-sz^\nu} z^{\nu-1}\,\mathrm{d}s .
\end{equation}
We now show that \eqref{iterg} satisfies equation \eqref{gtilde}; by
inserting the
Laplace transform into \eqref{gtilde}, we obtain
%
\begin{equation}
\label{insert-lapl}
z^\nu z^{\nu-1} \int_0^\infty G ( u,s )
\mathrm{e}^{-sz^\nu}\, \mathrm{d}s - z^{\nu-1} u
= ( \lambda u - \mu ) ( u-1 )
z^{\nu-1} \int_0^\infty\frac{\partial}{\partial u}
G ( u,s ) \mathrm{e}^{-sz^\nu} \,\mathrm{d}s .
\end{equation}
The inversion of the integral with $\partial/\partial u$ is justified
because
%
\begin{eqnarray}
\biggl| \frac{\partial}{\partial u} G (u,s) \biggr| &=&
\Biggl| \sum_{k=1}^\infty k u^{k-1} \Pr
\{ N(s) =k \} \Biggr| \nonumber\\[-8pt]\\[-8pt]
&\leq&\sum_{k=1}^\infty
k \Pr \{ N(s)=k \} = \mathbb{E} N(s) < \infty.\nonumber
\end{eqnarray}
Taking into account that $G ( u,t )$ satisfies the first-order
PDE
%
\begin{equation}
\frac{\partial G}{\partial s} =
( \lambda u -\mu ) ( u-1 ) \frac{\partial G}{
\partial u},
\end{equation}
from \eqref{insert-lapl}, we have that
%
\begin{eqnarray}
z^\nu\int_0^\infty G ( u,s ) \mathrm{e}^{-sz^\nu} \,\mathrm{d}s -u &=&
\int_0^\infty\frac{\partial}{\partial s} G ( u,s )
\mathrm{e}^{-sz^\nu} \,\mathrm{d}s \nonumber\\
& = & G ( u,s ) \mathrm{e}^{-sz^\nu} \bigl|_{s=0}^{s=\infty}
+ z^\nu\int_0^\infty G ( u,s ) \mathrm{e}^{-s z^\nu} \,\mathrm{d}s
\\
& = & -u + z^\nu\int_0^\infty G ( u,s ) \mathrm{e}^{-s z^\nu} \,\mathrm{d}s
\nonumber.
\end{eqnarray}
This shows that \eqref{ite} holds for the one-dimensional distributions.
This concludes the proof of Theorem \ref{teo-ite}.
\end{pf}

\begin{rem}
For $\nu=1/2^n, n \in\mathbb{N}$, the density $f_{T_{2 \nu}}$ of the
random time $T_{2 \nu}$ appearing in \eqref{ite} becomes the probability
density of an $ ( n-1 )$-iterated Brownian motion, that is,
\begin{eqnarray}
\label{ite-bm}
 \Pr \bigl\{ T_{{1}/{(2^{n-1})}} ( t ) \in \mathrm{d}s \bigr\} &=&
\Pr \{ \vert\mathcal{B}_1 ( \vert\mathcal{B}_2 (
\cdots\vert\mathcal{B}_n
( t ) \vert\cdots ) \vert
) \vert \in \mathrm{d}s \} \nonumber\\
& =& 2^n \int_0^\infty\frac{\mathrm{e}^{- s^2/{(4
\omega_1)}}}{\sqrt{4 \uppi\omega_1}} \,\mathrm{d} \omega_1
\int_0^\infty\frac{\mathrm{e}^{- \omega_1^2/{(4
\omega_2)}}}{\sqrt{4 \uppi\omega_2}} \,\mathrm{d} \omega_2 \cdots\\
&&{}\times
\int_0^\infty\frac{\mathrm{e}^{- \omega_{n-1}^2/{(4t)}}}{\sqrt{4
\uppi t}} \,\mathrm{d} \omega_{n-1} , \nonumber
\end{eqnarray}
as can easily be inferred from \citet{ors2008}, Theorem 2.1. The difference
between \eqref{ite-bm} and its corresponding formula in the cited
paper is that
here, the diffusion coefficient is equal to one.
\end{rem}

In the following theorems, we separately derive the three different
expressions of the
probability of extinction in the cases $\lambda>\mu, \lambda<\mu$ and
$\lambda=\mu$.
We prefer to treat them separately because their proofs are somewhat different.

\begin{thm}
For a fractional linear birth--death process with rates $\lambda>\mu$, the
probability of extinction has the form
%
\begin{eqnarray}
\label{exti}
p_0^\nu ( t ) &=& \Pr
\{ N_\nu ( t ) = 0 \} \nonumber\\[-8pt]\\[-8pt]
&=& \frac{\mu}{\lambda}
- \frac{\lambda- \mu}{\lambda}
\sum_{m=1}^{\infty} \biggl( \frac{\mu}{\lambda} \biggr)^m E_{\nu,1}
\bigl( -t^\nu ( \lambda- \mu) m \bigr)\nonumber
\end{eqnarray}
for $t>0, 0<\nu\leq1$, and where $E_{\nu,1} ( x )$ is the Mittag--Leffler
function \eqref{mittag}.
\end{thm}

\begin{pf}
In light of the subordination relationship \eqref{ite} of Theorem
\ref{teo-ite}, and by taking into account the extinction probability
of the classical linear birth--death process
%
\begin{equation}
\label{cla-ext}
\Pr \{ N ( t ) = 0 \} = \frac{\mu
- \mu \mathrm{e}^{-t ( \lambda
- \mu ) }}{\lambda- \mu \mathrm{e}^{-t ( \lambda- \mu ) }},\qquad
t>0,
\end{equation}
we can write that
%
\begin{equation}
\label{fract-ext-sub}
\Pr \{ N_\nu ( t ) = 0 \} = \int_0^{+ \infty}
\frac{\mu- \mu \mathrm{e}^{-s ( \lambda- \mu ) }}{\lambda- \mu \mathrm{e}^{-s
( \lambda- \mu ) }} \Pr \{T_{2 \nu}
( t ) \in \mathrm{d}s \}
\end{equation}
for all $t>0$ and $0 < \nu\leq1$. By taking the Laplace transform of
\eqref{fract-ext-sub}, we obtain that
\begin{eqnarray}\label{aba}
&& \int_0^\infty \mathrm{e}^{-zt} \Pr \{ N_\nu
( t ) = 0 \} \,\mathrm{d}t \nonumber\\
&&\quad  = \int_0^\infty\frac{\mu- \mu \mathrm{e}^{-s ( \lambda-\mu )}}{
\lambda-\mu \mathrm{e}^{-s ( \lambda- \mu )}} z^{\nu-1}
\mathrm{e}^{-sz^\nu} \,\mathrm{d}s \nonumber\\
&&\quad  =\frac{\mu}{\lambda} \int_0^\infty \bigl( 1-\mathrm{e}^{-s (\lambda
-\mu )} \bigr) \sum_{m=0}^\infty \biggl( \frac{\mu}{\lambda}
\mathrm{e}^{-s ( \lambda-\mu )} \biggr)^m z^{\nu-1} \mathrm{e}^{-sz^\nu} \,\mathrm{d}s
\nonumber\\[-8pt]\\[-8pt]
&&\quad  = \frac{\mu}{\lambda} \sum_{m=0}^\infty \biggl( \frac{\mu}{\lambda}
\biggr)^m z^{\nu-1} \biggl[ \int_0^\infty \bigl(
\mathrm{e}^{-s ( \lambda-\mu ) m -sz^\nu} - \mathrm{e}^{-s ( \lambda
- \mu ) ( 1+m ) -sz^\nu} \bigr) \,\mathrm{d}s \biggr]
\nonumber\\
&&\quad  = \frac{\mu}{\lambda} z^{\nu-1} \Biggl\{ \sum_{m=0}^\infty
\biggl( \frac{\mu}{\lambda} \biggr)^m \frac{1}{ ( \lambda-\mu
)m + z^\nu} - \sum_{m=1}^\infty \biggl( \frac{\mu}{\lambda}
\biggr)^{m-1} \frac{1}{ ( \lambda-\mu ) m +z^\nu}
\Biggr\} \nonumber\\
&&\quad  = \frac{\mu}{\lambda} z^{\nu-1} \sum_{m=1}^\infty
\frac{1}{ ( \lambda-\mu ) m +z^\nu} \biggl(
\frac{\mu}{\lambda} \biggr)^m \biggl( 1 - \frac{\lambda}{\mu} \biggr)
+ \frac{\mu}{\lambda} z^{\nu-1} \frac{1}{z^\nu} \nonumber.
\end{eqnarray}
The above steps are valid because $0<\frac{\mu}{\lambda}\mathrm{e}^{-s ( \lambda
-\mu
)}<1$ for $\lambda>\mu$.
By keeping in mind the Laplace transform of the Mittag--Leffler
function $E_{\nu,1} (-x t^\nu),$
%
\begin{equation}
\label{mitta-lap}
\int_0^\infty \mathrm{e}^{-st} E_{\nu,1} (-x t^\nu) \,\mathrm{d}t = \frac{s^{\nu-1}}{
s^\nu+x},
\end{equation}
we readily arrive at the claimed result.
\end{pf}
\begin{rem}
When $\nu=1$, we obtain from \eqref{exti} the form of
the extinction probability \eqref{cla-ext}
for the classical birth--death model:
%
\begin{eqnarray}
\label{eq:classical-ext}
 \Pr \{ N ( t ) = 0 \} &=& \frac{\mu
- \lambda}{\lambda} \Biggl[
\sum_{m=1}^{+ \infty} \biggl( \frac{\mu}{\lambda} \biggr)^m \mathrm{e}^{- (
\lambda- \mu )
m t} \Biggr] + \frac{\mu}{\lambda}\nonumber\\
& = & \frac{\mu- \lambda}{\lambda} \biggl[ \frac{1}{1 -
(\mu/{\lambda}) \mathrm{e}^{- (
\lambda- \mu ) t }} - 1 \biggr] + \frac{\mu}{\lambda}
\nonumber\\[-8pt]\\[-8pt]
&=&
\frac{\mu- \lambda}{\lambda} \biggl[ \frac{({\mu}/{\lambda})
\mathrm{e}^{-t ( \lambda- \mu )}}{
1- ({\mu}/{\lambda}) \mathrm{e}^{-t ( \lambda- \mu )}}
\biggr] + \frac{\mu}{\lambda}\nonumber\\
&=& \frac{\mu- \mu \mathrm{e}^{-t ( \lambda- \mu )}}{ \lambda
- \mu \mathrm{e}^{-t ( \lambda- \mu )}}. \nonumber
\end{eqnarray}
From \eqref{fract-ext-sub} for $\nu=1$, $\Pr \{ T_2 ( t )
\in \mathrm{d}s \} = \delta ( s - t )$ and we again retrieve result
\eqref{cla-ext}.
\end{rem}

\begin{rem}
From \eqref{exti}, we note that
%
\begin{equation}
\label{eq:limit-behaviour}
\Pr \{ N_\nu ( t ) = 0 \}
\stackrel{t \rightarrow+ \infty}{\longrightarrow}
\frac{\mu}{\lambda}\qquad  \forall \nu\in (0, 1 ],
\end{equation}
which is the asymptotic extinction probability, irrespective of the
value of $\nu$.
\end{rem}

Let us now deal with the case $\lambda< \mu$, that is, when the rate
of birth is strictly
lower than the rate of death.

\begin{thm}
For $\mu>\lambda$, the probability $p_0^\nu ( t ) = \Pr
\{ N_\nu (
t ) = 0 \}$ of
complete extinction of the population is
%
\begin{equation}
\label{eq:prob-ext-mu-great}
p_0^\nu ( t ) = 1 - \frac{\mu- \lambda}{\lambda}
\sum_{m=1}^{+\infty} \biggl( \frac{\lambda}{\mu} \biggr)^m E_{\nu,1}
\bigl( - t^\nu ( \mu- \lambda ) m \bigr),
\end{equation}
where $t>0$, $0 < \nu\leq1$ and $E_{\nu,1} ( x )$ is the Mittag--Leffler
function \eqref{mittag}.
\end{thm}

\begin{pf}
We start by rewriting \eqref{cla-ext} as
%
\begin{equation}
\label{eq:new-clas}
p_0 ( t ) = \dfrac{\mu \mathrm{e}^{-t ( \mu- \lambda )}
- \mu}{\lambda \mathrm{e}^{-t ( \mu- \lambda )} - \mu} .
\end{equation}
Using \eqref{ite}, we are able to write
%
\begin{equation}
\label{eq:mu-great}
p_0^\nu ( t ) = \int_0^{+ \infty} \frac{\mu \mathrm{e}^{-s (
\mu- \lambda ) } - \mu}{\lambda \mathrm{e}^{-s ( \mu- \lambda )}
- \mu} \Pr \{ T_{2 \nu} ( t ) \in \mathrm{d}s \} .
\end{equation}
By applying the Laplace transform to \eqref{eq:mu-great}, we obtain that
\begin{eqnarray}
\label{laplace-tr}
L_0^\nu ( z ) & =& \int_0^{+ \infty} \frac{ \mu \mathrm{e}^{-s (
\mu- \lambda )} - \mu}{\lambda \mathrm{e}^{-s ( \mu- \lambda )}
- \mu} z^{\nu-1} \mathrm{e}^{-sz^\nu} \,\mathrm{d}s\nonumber\\
& =&
\int_0^{+ \infty} \frac{\mathrm{e}^{-s ( \mu- \lambda )} -1}{
({\lambda}/{\mu}) \mathrm{e}^{-s ( \mu- \lambda )} -1}
z^{\nu-1} \mathrm{e}^{-s z^\nu} \,\mathrm{d}s \nonumber\\
& =& z^{\nu-1} \int_0^{+ \infty} \bigl( 1-\mathrm{e}^{- s
( \mu- \lambda )}
\bigr) \mathrm{e}^{-s z^\nu} \sum_{m=0}^{+ \infty} \biggl[ \frac{\lambda}{\mu} \mathrm{e}^{-s
( \mu- \lambda )} \biggr]^m
\nonumber\\
& =& z^{\nu-1} \sum_{m=0}^{+ \infty} \biggl( \frac{\lambda}{\mu} \biggr)^m
\int_0^{+ \infty} \bigl( 1-\mathrm{e}^{-s ( \mu- \lambda ) }
\bigr) \mathrm{e}^{-sz^\nu} \mathrm{e}^{-s ( \mu- \lambda ) m} \,\mathrm{d}s
\nonumber\\
& = &z^{\nu-1} \sum_{m=0}^{+ \infty} \biggl( \frac{\lambda}{\mu} \biggr)^m
\int_0^{+ \infty} \mathrm{e}^{-s ( \mu- \lambda ) m -sz^\nu} -
\mathrm{e}^{-s ( \mu- \lambda ) ( m+1 ) -sz^\nu} \mathrm{d}s
\nonumber\\[-8pt]\\[-8pt]
& =& z^{\nu-1} \sum_{m=0}^{+\infty} \biggl( \frac{\lambda}{\mu} \biggr)^m
\biggl\{ \frac{1}{ ( \mu- \lambda ) m + z^\nu} -
\frac{1}{ ( \mu- \lambda ) ( m+1 )
+ z^\nu} \biggr\} \nonumber\\
& =& z^{\nu-1} \Biggl\{ \sum_{m=0}^{+ \infty} \biggl( \frac{\lambda}{\mu}
\biggr)^m \frac{1}{ ( \mu- \lambda ) m+z^\nu} -
\sum_{m=1}^{+ \infty} \biggl( \frac{\lambda}{\mu} \biggr)^{m-1}
\frac{1}{ ( \mu- \lambda ) m + z^\nu} \Biggr\}
\nonumber\\
& =& z^{\nu-1} \Biggl\{ \frac{1}{z^\nu} + \sum_{m=1}^{+ \infty}
\biggl( \frac{\lambda}{\mu} \biggr)^m \frac{1}{ ( \mu- \lambda
) m + z^\nu} - \frac{\mu}{\lambda} \sum_{m=1}^{+ \infty}
\biggl( \frac{\lambda}{\mu} \biggr)^m \frac{1}{ ( \mu- \lambda
) m +z^\nu} \Biggr\} \nonumber\\
& =& z^{\nu-1} \Biggl\{ \frac{1}{z^\nu} + \biggl[ 1-\frac{\mu}{\lambda}
\biggr] \sum_{m=1}^{+ \infty} \biggl( \frac{\lambda}{\mu} \biggr)^m
\frac{1}{ ( \mu- \lambda ) m +z^\nu} \Biggr\} \nonumber\\
& =& \frac{1}{z} + \biggl[ 1- \frac{\mu}{\lambda} \biggr]
\sum_{m=1}^{+ \infty} \biggl( \frac{\lambda}{\mu} \biggr)^m
\frac{z^{\nu-1}}{ ( \mu- \lambda ) m + z^\nu} \nonumber.
\end{eqnarray}
Inverting \eqref{laplace-tr} by means of \eqref{mitta-lap}, we retrieve formula
\eqref{eq:prob-ext-mu-great}.
\end{pf}

\begin{rem}
When $\nu=1$, we reobtain from \eqref{eq:prob-ext-mu-great} the extinction
probability of the classical birth--death process \eqref{eq:new-clas}:
%
\begin{eqnarray}\label{aa}
p_0^1 ( t ) & = & 1 - \biggl[ \frac{\mu- \lambda}{\lambda} \biggr]
\sum_{m=1}^{+ \infty} \biggl( \frac{ \lambda}{\mu} \biggr)^m \mathrm{e}^{- (
\mu- \lambda ) m t}\nonumber \\
& = & 1- \biggl( \frac{\mu- \lambda}{\lambda} \biggr) \biggl( \frac{1}{
1-(\lambda/\mu) \mathrm{e}^{- ( \mu- \lambda ) t}} -1 \biggr)
\\
& = & 1 - \biggl( \frac{\mu- \lambda}{\lambda} \biggr) \biggl[
\frac{(\lambda/\mu) \mathrm{e}^{- \lambda ( \mu- \lambda ) t}}{
1 - (\lambda/\mu) \mathrm{e}^{- ( \mu- \lambda ) t} } \biggr]\nonumber
\\
& = & 1+ \frac{(\lambda^2/\mu) \mathrm{e}^{- ( \mu- \lambda ) t} -
\lambda \mathrm{e}^{- ( \mu- \lambda ) t}}{\lambda- (\lambda^2/\mu)
\mathrm{e}^{- ( \mu- \lambda ) t} } =
\frac{\lambda- \lambda \mathrm{e}^{- ( \mu- \lambda ) t}}{
\lambda- (\lambda^2/\mu) \mathrm{e}^{- ( \mu- \lambda ) t}}
\nonumber\\
& = & \frac{1 - \mathrm{e}^{- ( \mu- \lambda ) t}}{1 - (\lambda/\mu)
\mathrm{e}^{- ( \mu- \lambda ) t}} =
\frac{\mu \mathrm{e}^{-t ( \mu- \lambda ) } - \mu}{\lambda \mathrm{e}^{-t
( \mu- \lambda ) } - \mu} \nonumber.
\end{eqnarray}
\end{rem}

\begin{rem}
Population extinction in the long run is evident from \eqref
{eq:prob-ext-mu-great}
as
%
\begin{equation}
\Pr \{ N_\nu ( t ) = 0 \}
\stackrel{t \rightarrow+ \infty}{\longrightarrow} 1,
\end{equation}
due to the death rate exceeding the birth rate for all
$0<\nu\leq1$.
\end{rem}

In the next theorem, we treat the remaining case, that is, when $\mu
=\lambda$.

\begin{thm}
For the fractional linear birth process, when the rate of birth equals
the rate
of death (i.e., when $\lambda=\mu$), the extinction probability
$p_0^\nu ( t )$
reads
%
\begin{equation}
\label{eq:ext-prob-equal}
p_0^\nu ( t ) =
\frac{\lambda t^\nu}{\nu} \int_0^{+ \infty}
\mathrm{e}^{-w} E_{\nu,\nu} ( -w \lambda t^\nu ) \,\mathrm{d}w
= 1 - \int_0^{+ \infty} \mathrm{e}^{-w} E_{\nu,1} (
- \lambda t^\nu w ) \,\mathrm{d}w
\end{equation}
with $t>0$, $0 < \nu\leq1$ and where $E_{\nu,1} ( x )$ is the Mittag--Leffler
function \eqref{mittag}.
\end{thm}

\begin{pf}
Again using \eqref{ite}, we write
%
\begin{equation}
\label{eq:equal}
p_0^\nu ( t ) = \int_0^{+ \infty} \frac{\lambda s}{1+
\lambda s} \Pr \{ T_{2 \nu} ( t ) \in \mathrm{d}s \} .
\end{equation}
We now apply the Laplace transform once again, thus obtaining
%
\begin{eqnarray}
L_0^\nu ( z ) & =& \int_0^{+ \infty} \frac{ \lambda s z^{\nu-1}
\mathrm{e}^{- z^\nu s}}{\lambda s +1} \,\mathrm{d}s \nonumber\\
& =& \lambda z^{\nu-1} \int_0^{+ \infty} s \mathrm{e}^{- z^\nu s} \int_0^{+
\infty}
\mathrm{e}^{- w ( \lambda s +1 )} \,\mathrm{d}w \,\mathrm{d}s \nonumber\\
& = &\lambda z^{\nu-1} \int_0^{+ \infty} \mathrm{e}^{-w} \int_0^{+ \infty}
s \mathrm{e}^{- z^\nu s - w \lambda s} \,\mathrm{d}s \,\mathrm{d}w
\\
&
\stackrel{ ( y=s ( z^\nu+ \lambda w ) )}{=}&
\lambda z^{\nu-1} \int_0^{+ \infty} \mathrm{e}^{-w} \int_0^{+ \infty}
\frac{y}{ z^\nu+ \lambda w } \mathrm{e}^{-y} \frac{\mathrm{d}y}{
z^\nu+ \lambda w } \,\mathrm{d}w \nonumber\\
& =& \lambda\int_0^{+ \infty} \mathrm{e}^{-w} \frac{1}{z^\nu+ \lambda w}
\cdot\frac{z^{\nu-1}}{z^\nu+ \lambda w}\, \mathrm{d}w \nonumber.
\end{eqnarray}
By inverting the Laplace transform, we obtain the integral
form
%
\begin{equation}
\label{eq:ext-prob-int}
p_0^\nu ( t ) = \lambda\int_0^{+ \infty} \mathrm{e}^{- w} \int_0^t
u^{\nu-1} E_{\nu,\nu} ( -w \lambda u^\nu ) E_{\nu,1} (
- w \lambda ( t-u )^\nu ) \,\mathrm{d}u \,\mathrm{d}w ,
\end{equation}
which involves convolutions of generalized
Mittag--Leffler functions $E_{\alpha, \beta} ( t )$,
defined, for example, in \citet{podlubny}, equation (1.56), page 17.
The inner integral in \eqref{eq:ext-prob-int} can be worked out explicitly
as follows:
%
\begin{eqnarray}\label{conv}
&& \int_0^t u^{\nu-1} E_{\nu,\nu} ( -w \lambda u^\nu )
E_{\nu,1} \bigl( -w \lambda ( t-u )^\nu \bigr) \,\mathrm{d}u \nonumber\\
&&\qquad\hspace*{2pt}  = \sum_{m=0}^{\infty} \sum_{r=0}^{\infty} \frac{ ( -w \lambda
)^m}{\Gamma ( \nu m + \nu )} \frac{
(- w \lambda )^r}{\Gamma ( \nu r +1 )}
\int_0^t u^{\nu-1} u^{\nu m} ( t-u )^{\nu r} \,\mathrm{d}u
\nonumber\\
&&\qquad\hspace*{2pt}  = \sum_{m=0}^{\infty} \sum_{r=0}^{\infty} \frac{ ( -w \lambda
)^{m+r}}{\Gamma ( \nu m + \nu ) \Gamma (
\nu r +1 )} t^{\nu+ \nu ( m+r )} \frac{
\Gamma ( \nu m + \nu ) \Gamma ( \nu r +1 )}{
\Gamma ( \nu ( m+r ) + \nu+1 )}
\\
&&\quad  \stackrel{(m+r=n)}{=}
\sum_{m=0}^\infty\sum_{n=m}^\infty\frac{ ( -w \lambda
)^n}{\Gamma ( \nu n + \nu+1 )} t^{\nu+ \nu n}
= \sum_{n=0}^\infty\frac{ ( -w \lambda )^n}{
\Gamma ( \nu n + \nu+1 )} t^{\nu+ \nu n} (
n+1 ) \nonumber\\
&&\qquad\hspace*{2pt}  = \frac{t^\nu}{\nu} \sum_{n=0}^{\infty}
\frac{ (- w \lambda t^\nu )^n}{\Gamma ( \nu
( n+1 ) )} =
\frac{t^\nu}{\nu} E_{\nu,\nu} ( -w \lambda t^\nu ) .
\nonumber
\end{eqnarray}
The extinction probability now reads
%
\begin{equation}
\label{first-form}
p_0^\nu ( t ) = \frac{\lambda t^\nu}{\nu} \int_0^\infty
\mathrm{e}^{-w} E_{\nu,\nu} ( -w \lambda t^\nu ) \,\mathrm{d}w .
\end{equation}
Using the relationship
%
\begin{equation}
\label{relmitt}
\frac{\mathrm{d}}{\mathrm{d}x} E_{\nu,1} ( x ) = \frac{1}{\nu} E_{\nu,\nu}
( x ),
\end{equation}
the extinction probability \eqref{first-form} takes the alternative form
\eqref{eq:ext-prob-equal} because
%
\begin{eqnarray}
 p_0^\nu ( t ) &\stackrel{ (
-w \lambda t^\nu=y )}{=}&
- \frac{\lambda t^\nu}{\nu} \int_0^{- \infty} E_{\nu,\nu} ( y )
\mathrm{e}^{y/({\lambda t^\nu}) } \,\mathrm{d}y \nonumber\\
& =& \frac{1}{\nu} \int_{-\infty}^{0} E_{\nu,\nu} ( y
) \mathrm{e}^{{y}/{(\lambda t^\nu)}} \,\mathrm{d}y \stackrel{\mathrm{\eqref{relmitt}}}{
=} \int_{-\infty}^0 \mathrm{e}^{y/{(\lambda t^\nu)}} \frac{\mathrm{d}}{\mathrm{d}y} E_{
\nu,1} ( y ) \,\mathrm{d}y \\
& =& 1- \frac{1}{\lambda t^\nu} \int_{-\infty}^0
\mathrm{e}^{y/{(\lambda t^\nu)}}
E_{\nu,1} ( y ) \,\mathrm{d}y \stackrel{ ( w=- y/{(\lambda
t^\nu)} )}{=}
1- \int_0^\infty \mathrm{e}^{-w} E_{\nu,1} ( -\lambda t^\nu w ) \,\mathrm{d}w.
\nonumber
\end{eqnarray}
This completes the proof of \eqref{eq:ext-prob-equal}.
\end{pf}

\begin{rem}
From \eqref{eq:ext-prob-equal}, when $\nu=1$, we again retrieve the
classical form
%
\begin{equation}
p_0 ( t ) = \frac{\lambda t}{\lambda t +1},
\end{equation}
as expected.
\end{rem}
\begin{rem}
The limiting extinction probability when $\mu=\lambda$ is
%
\begin{equation}
\Pr \{ N_\nu ( t ) = 0 \}
\stackrel{t \rightarrow+ \infty}{\longrightarrow} 1
\end{equation}
for all values of $0<\nu\leq1$.
\end{rem}

\begin{rem}
The last expression in \eqref{eq:ext-prob-equal}
is in some ways similar to the {Riemann} limit for $\mu\rightarrow
\lambda$ of \eqref{exti}
and \eqref{eq:prob-ext-mu-great}.
\end{rem}
\begin{rem}
We can rewrite the probabilities \eqref{total-extinction} in an alternative
form which permits us to give an interesting interpretation to their structure.

For the case $\lambda>\mu$, we can write
\begin{eqnarray}
p_0^\nu ( t ) & =& \frac{\mu}{\lambda} \Biggl[
1-\frac{\lambda}{\mu} \frac{\lambda-\mu}{\lambda} \sum_{m=1}^\infty
\biggl( \frac{\mu}{\lambda} \biggr)^m E_{\nu,1} \bigl( - t^\nu
( \lambda-\mu ) m \bigr) \Biggr] \nonumber\\[-8pt]\\[-8pt]
& =& \frac{\mu}{\lambda} \Biggl[ 1- \sum_{m=1}^\infty\Pr
\{ \mathcal{G}=m | \mathcal{G} \geq1 \} E_{\nu,1}
\bigl( -t^\nu ( \lambda-\mu ) m \bigr) \Biggr], \nonumber
\end{eqnarray}
where $\mathcal{G}$ is a geometric r.v. with distribution
%
\begin{equation}
\Pr ( \mathcal{G}=m | \mathcal{G} \geq1 ) =
\frac{\Pr ( \mathcal{G} = m )}{\Pr
( \mathcal{G} \geq1 )} =
\frac{\lambda-\mu}{\lambda} \biggl( \frac{\mu}{\lambda} \biggr)^m
\frac{\lambda}{\mu},\qquad  m \geq1 .
\end{equation}

The treatment of the opposite case $\lambda<\mu$ is similar except
that a different
conditional geometric r.v. $\mathcal{G}'$ must be introduced, defined as
%
\begin{equation}
\Pr ( \mathcal{G}'=m | \mathcal{G}' \geq1 ) =
\frac{\mu}{\lambda} \biggl( \frac{\lambda}{\mu} \biggr)^m
\frac{\mu-\lambda}{\mu},\qquad  m \geq1,
\end{equation}
and thus
%
\begin{equation}
p_0^\nu ( t ) = 1- \sum_{m=1}^\infty\Pr
( \mathcal{G}'=m | \mathcal{G}' \geq1 ) E_{\nu,1} \bigl(
-t^\nu ( \mu-\lambda ) m \bigr) .
\end{equation}

A well-known property for a fractional Poisson process $\mathcal
{N}_\nu (
t ), t > 0$,
of degree $0<\nu\leq1$ and parameter $\lambda>0$ is that \citet{fr-poisson}
%
\begin{equation}
\Pr \{ \mathcal{N}_\nu ( t ) =0 \} =
E_{\nu,1} ( -t^\nu\lambda ) = \Pr
( \mathcal{T}_\nu\geq t ),
\end{equation}
where $\mathcal{T}_\nu= \inf ( s \dvtx\mathcal{N}_\nu ( s ) =1
)$. This permits us to rewrite the extinction probabilities also in terms
of waiting times of a fractional Poisson process with a random rate
$\lambda\mathcal{G}$.

For the case $\lambda=\mu$, the interpretation is straightforward
because the waiting time of the related fractional Poisson process
has a rate $\lambda\mathcal{E}$, where $\mathcal{E}$ is an exponentially
distributed r.v. with parameter equal to one.
\end{rem}

\begin{rem}
In the case $\mu=\lambda$, it is well known that the extinction
probability in the classical birth--death
process, $p_0 ( s )$, $s>0$, satisfies the nonlinear Riccati differential
equation
%
\begin{equation}
\label{riccati}
p_0' ( s ) + 2 \lambda p_0 ( s ) = \lambda+
\lambda [ p_0 ( s ) ]^2 .
\end{equation}
By using \eqref{riccati}, we can provide an alternative proof for the
subordination
relationship \eqref{ite}:
%
\begin{eqnarray}
&& \int_0^\infty p_0' ( s ) \Pr \bigl( T_{2 \nu} ( t )
\in \mathrm{d}s \bigr) + 2 \lambda p_0^\nu ( t ) =
\lambda+ \lambda\int_0^\infty [ p_0 ( s ) ]^2
\Pr \bigl( T_{2 \nu} ( t )
\in \mathrm{d}s \bigr) \\
&&\quad  \Longleftrightarrow\quad \int_0^\infty\frac{\lambda}{ ( 1 + \lambda s )^2 }
\Pr \bigl( T_{2 \nu} ( t )
\in \mathrm{d}s \bigr) + 2 \lambda p_0^\nu ( t ) =
\lambda+ \lambda\int_0^\infty\frac{\lambda^2 s^2}{ (
1+\lambda s )^2} \Pr \bigl( T_{2 \nu} ( t )
\in \mathrm{d}s \bigr) \nonumber\\
&&\quad  \Longleftrightarrow\quad \int_0^\infty\lambda\frac{ ( 1-\lambda^2 s^2 )}{
( 1+ \lambda s )^2} \Pr \bigl( T_{2 \nu} ( t
) \in \mathrm{d}s \bigr) = \lambda- 2 \lambda p_0^\nu ( t )
\nonumber\\
&&\quad  \Longleftrightarrow\quad \int_0^\infty\frac{1- \lambda s}{1+ \lambda s}
\Pr \bigl( T_{2 \nu} ( t ) \in \mathrm{d}s \bigr)
= 1 - 2 p_0^\nu ( t ) \nonumber\\
&&\quad  \Longleftrightarrow\quad 2 p_0^\nu ( t ) = 2 \int_0^\infty\frac{\lambda s}{
1+\lambda s} \Pr \bigl( T_{2 \nu} ( t ) \in \mathrm{d}s \bigr)
\nonumber\\
&&\quad  \Longleftrightarrow\quad  p_0^\nu ( t ) = \int_0^\infty\frac{\lambda s}{
1+\lambda s} \Pr \bigl( T_{2 \nu} ( t ) \in \mathrm{d}s \bigr) .
\end{eqnarray}
\end{rem}

\begin{rem}
By exploiting the subordination relationship \eqref{ite} and the fact
that the
extinction probability in the classical case satisfies the integral equation
%
\begin{equation}
\label{integral-eq-clas}
p_0 ( t ) = \int_0^t \mathrm{e}^{- ( \lambda+ \mu ) u} \mu \,\mathrm{d}u +
\int_0^t \lambda \mathrm{e}^{- ( \lambda+ \mu ) u} [ p_0 (
t-u ) ]^2 \,\mathrm{d}u,
\end{equation}
we can give an integral form for $p_0^\nu ( t )$:
%
\begin{equation}
\label{intform}
p_0^\nu ( t ) = \int_0^{+ \infty} \biggl\{ \int_0^s
\mathrm{e}^{- ( \lambda+ \mu ) u} \mu \,\mathrm{d}u + \int_0^s \lambda \mathrm{e}^{- (
\lambda+ \mu ) u} [p_0 ( s-u ) ]^2
\,\mathrm{d}u \biggr\} \Pr
\{ T_{2 \nu} ( t ) \in \mathrm{d}s \} .
\end{equation}
We note that the first integral of \eqref{intform} can be worked out
explicitly as
follows:
%
\begin{equation}
\mu\int_0^\infty \mathrm{e}^{-zt} \biggl[ \int_0^\infty\int_0^s
\mathrm{e}^{- ( \lambda+\mu ) u} \Pr \{ T_{2 \nu} ( t
) \in \mathrm{d}s \} \,\mathrm{d}u \biggr] \,\mathrm{d}t
= \frac{\mu}{z} \frac{1}{\lambda+ \mu+ z^\nu} .
\end{equation}
This can be directly inverted so as to obtain
%
\begin{eqnarray}
&& \int_0^\infty\int_0^s \mathrm{e}^{- ( \lambda+ \mu ) u} \mu
\Pr \{T_{2 \nu} ( t ) \in \mathrm{d}s \} \,\mathrm{d}u\nonumber\\
 &&\quad =
\mu\int_0^t w^{\nu-1} E_{\nu,\nu} \bigl( - ( \lambda+\mu ) w^\nu
\bigr) \,\mathrm{d}w \nonumber\\[-8pt]\\[-8pt]
&&\quad  =  \frac{\mu t^\nu}{\nu} \sum_{m=0}^\infty\frac{ (
- ( \lambda+ \mu ) t^\nu )^m}{ ( m+1 )
\Gamma ( \nu m+\nu )}\nonumber\\
&&\quad = \frac{\mu}{\lambda+\mu} \bigl[ 1- E_{\nu,1} \bigl(- (\lambda+ \mu)
t^\nu \bigr) \bigr] . \nonumber
\end{eqnarray}
\end{rem}

\section{The state probabilities of the fractional linear birth--death process}\label{sec2}

Here, we present three theorems concerning the structure of
the state probabilities
$\Pr \{ N_\nu ( t ) = k \}$, $t > 0,$
with $0 < \nu\leq1$. Three cases must be distinguished and treated
separately, as
in Section \ref{sec}, namely $\lambda>\mu$, $\lambda<\mu$ and $\lambda
=\mu$.
\begin{thm}
For the case $\lambda>\mu$,
the state probabilities $p_k^\nu ( t )$, $k \geq1$,
$t>0$, $0< \nu\leq1$, in the
fractional linear birth--death process $N_\nu ( t )$, $t>0$, have the
following form:
%
\begin{eqnarray}
\label{pklgm}
p_k^\nu ( t ) &=& \biggl( \frac{\lambda-\mu}{\lambda} \biggr)^2
\sum_{l=0}^\infty\pmatrix{l+k\cr l} \biggl( \frac{\mu}{\lambda} \biggr)^l
\sum_{r=0}^{k-1} ( -1 )^r \pmatrix{k-1\cr r}\nonumber\\[-8pt]\\[-8pt]
&&\hphantom{\biggl( \frac{\lambda-\mu}{\lambda} \biggr)^2
\sum_{l=0}^\infty\pmatrix{l+k\cr l} \biggl( \frac{\mu}{\lambda} \biggr)^l
\sum_{r=0}^{k-1}}
{}\times E_{\nu,1}\bigl ( - ( l+r+1 ) ( \lambda-\mu ) t^\nu
\bigr).\nonumber
\end{eqnarray}
\end{thm}

\begin{pf}
By exploiting the subordination
relationship \eqref{ite} and conveniently rewriting the
well-known form of the state probabilities of the classical linear birth--death
process, we have that
%
\begin{equation}
p_k^\nu ( t ) =
( \lambda-\mu )^2
\lambda^{k-1}
\int_0^\infty
\mathrm{e}^{- ( \lambda-\mu ) s}
\frac{ ( 1-\mathrm{e}^{- ( \lambda-\mu ) s}
)^{k-1}}{ ( \lambda-\mu \mathrm{e}^{- ( \lambda-\mu ) s}
)^{k+1}} \Pr \bigl( T_{2 \nu} (t) \in \mathrm{d}s \bigr) .
\end{equation}
By applying the Laplace transform, we obtain
%
\begin{eqnarray}
L_k^\nu ( z )& = & ( \lambda- \mu )^2
\lambda^{k-1} \int_0^\infty \mathrm{e}^{- ( \lambda-\mu ) s}
\frac{ ( 1-\mathrm{e}^{- ( \lambda-\mu ) s} )^{k-1}}{
( \lambda-\mu \mathrm{e}^{- ( \lambda- \mu ) s} )^{k+1}}
z^{\nu-1} \mathrm{e}^{-sz^\nu} \,\mathrm{d}s\nonumber \\
&= & ( \lambda-\mu )^2 \lambda^{k-1}
\sum_{l=0}^\infty\sum_{r=0}^{k-1} \pmatrix{- (k+1 )\cr l}
(-1 )^l \biggl(\frac{\mu}{\lambda} \biggr)^l \lambda^{-
( k+1 )} \pmatrix{k-1\cr r} ( -1 )^r
z^{\nu-1}\nonumber\\
&&\hphantom{( \lambda-\mu )^2 \lambda^{k-1}
\sum_{l=0}^\infty\sum_{r=0}^{k-1}}{} \times \int_0^\infty \mathrm{e}^{-sz^\nu} \mathrm{e}^{- (\lambda-\mu )sl}
\mathrm{e}^{- (\lambda-\mu ) sr} \mathrm{e}^{- ( \lambda-\mu ) s}
\,\mathrm{d}s \\
&=  & \biggl( \frac{\lambda-\mu}{\lambda} \biggr)^2 \sum_{l=0}^{\infty}
\sum_{r=0}^{k-1} \pmatrix{l+k\cr l} \pmatrix{k-1\cr r} (-1 )^r
\biggl( \frac{\mu}{\lambda} \biggr)^l z^{\nu-1} \int_0^\infty
\mathrm{e}^{-s ( z^\nu+ ( \lambda-\mu ) ( l+r+1 )
)} \,\mathrm{d}s \nonumber\\
& =  & \biggl( \frac{\lambda-\mu}{\lambda} \biggr)^2
\sum_{l=0}^\infty\sum_{r=0}^{k-1} \pmatrix{l+k\cr l} \pmatrix{k-1\cr r}
(-1 )^r \biggl( \frac{\mu}{\lambda} \biggr)^l
\frac{z^{\nu-1}}{z^\nu+ (\lambda- \mu ) (
l+r+1 )}, \nonumber
\end{eqnarray}
which can be easily inverted by using \eqref{mitta-lap},
thus obtaining \eqref{pklgm}.
\end{pf}

\begin{rem}
We check that for $\nu=1$, formula \eqref{pklgm} converts into the
well-known distribution of the linear birth--death process, thus being its
fractional extension.
For $\nu=1$, we get from \eqref{pklgm} that
%
\begin{equation}
p_k^1 ( t ) = \biggl( \frac{\lambda- \mu}{\lambda} \biggr)^2
\sum_{l=0}^\infty\pmatrix{l+k\cr l} \biggl( \frac{\mu}{\lambda} \biggr)^l
\sum_{r=0}^{k-1} ( -1 )^r \pmatrix{k-1\cr r}
\mathrm{e}^{- ( \lambda- \mu ) t ( l+r+1 )} .
\end{equation}
We now observe that
%
\begin{eqnarray}
\label{great1}
\sum_{r=0}^{k-1} ( -1 )^r \pmatrix{k-1\cr r} \mathrm{e}^{-t (
\lambda- \mu ) r} &=& \bigl( 1- \mathrm{e}^{- ( \lambda- \mu ) t}
\bigr)^{k-1},
\\
\label{great}
\sum_{l=0}^\infty\pmatrix{l+k\cr l} \biggl( \frac{\mu}{\lambda}
\biggr)^l \mathrm{e}^{- ( \lambda- \mu ) tl} &=& \biggl(
1- \frac{\mu}{\lambda} \mathrm{e}^{- ( \lambda- \mu ) t} \biggr)^{-
( k+1 ) },
\end{eqnarray}
where, in \eqref{great}, we applied the binomial expression
%
\begin{equation}
\sum_{l=0}^\infty\pmatrix{a+l\cr l} b^l =
\sum_{l=0}^\infty\pmatrix{- ( a+1 )\cr l} ( -b )^l
= ( 1-b )^{- ( a+1 )} .
\end{equation}
This permits us to write
%
\begin{equation}
p_k^1 ( t ) = \biggl( \frac{\lambda- \mu}{\lambda} \biggr)^2
\mathrm{e}^{- ( \lambda- \mu ) t} \frac{ ( 1 - \mathrm{e}^{- (
\lambda- \mu ) t} )^{k-1}}{ ( 1- (\mu/{\lambda})
\mathrm{e}^{- ( \lambda- \mu ) t} )^{k+1}},\qquad  \mu< \lambda,
\end{equation}
which coincides with \eqref{class-state}.
\end{rem}

\begin{rem}
In order to prove that $\sum_{k=0}^\infty p_k^\nu(t) = 1$ for $\lambda
>\mu$ (formula
\eqref{pklgm}), we can again apply the Laplace transform and prove that
$ \sum_{k=0}^\infty\int_0^\infty \mathrm{e}^{-zt} p_k^\nu(t) \,\mathrm{d}t = 1/z$. We
first calculate
\begin{eqnarray}
&& \sum_{k=1}^\infty\int_0^\infty \mathrm{e}^{-z t} p_k^\nu(t) \,\mathrm{d}t \nonumber\\
&&\quad  = \sum_{k=1}^\infty \biggl( \frac{\lambda-\mu}{\lambda} \biggr)^2
\sum_{l=0}^\infty\sum_{r=0}^{k-1} \pmatrix{l+k\cr l} \pmatrix{k-1\cr r}
(-1 )^r \biggl( \frac{\mu}{\lambda} \biggr)^l
\frac{z^{\nu-1}}{z^\nu+ (\lambda- \mu ) (
l+r+1 )} \qquad \qquad \\ 
&&\quad  = \sum_{k=1}^\infty \biggl( \frac{\lambda-\mu}{\lambda} \biggr)^2
\sum_{l=0}^\infty\sum_{r=0}^{k-1} \pmatrix{- ( k+1 )\cr l}
(-1)^l \pmatrix{k-1\cr r} (-1)^r \biggl( \frac{\mu}{\lambda} \biggr)^l
z^{\nu-1}\nonumber\\
&&\qquad \hphantom{\sum_{k=1}^\infty \biggl( \frac{\lambda-\mu}{\lambda} \biggr)^2\sum_{l=0}^\infty\sum_{r=0}^{k-1}}
{}\times \int_0^\infty \mathrm{e}^{-sz^\nu} \mathrm{e}^{-ls ( \lambda-\mu
)} \mathrm{e}^{-sr ( \lambda-\mu )} \mathrm{e}^{-s ( \lambda
- \mu )} \,\mathrm{d}s . \nonumber
\end{eqnarray}
By keeping in mind formulas \eqref{great1} and \eqref{great}, we have that
%
\begin{eqnarray}
&&\sum_{k=1}^\infty\int_0^\infty \mathrm{e}^{-z t} p_k^\nu(t)
\,\mathrm{d}t\nonumber\\
&&\quad  = \sum_{k=1}^\infty(\lambda-\mu)^2 \lambda^{k-1} \int_0^\infty
\mathrm{e}^{-s ( \lambda-\mu )} \frac{ ( 1-\mathrm{e}^{-s ( \lambda
-\mu )} )^{k-1} }{ ( \lambda-\mu \mathrm{e}^{-s ( \lambda- \mu
)} )^{k+1}} z^{\nu-1} \mathrm{e}^{-sz^\nu}\, \mathrm{d}s \\
&&\quad  = ( \lambda-\mu ) z^{\nu-1} \int_0^\infty
\frac{\mathrm{e}^{-sz^\nu}}{\lambda- \mu \mathrm{e}^{-s ( \lambda-\mu )}} . \nonumber
\end{eqnarray}
By using the Laplace transform of the extinction probability
(second line of formula \eqref{aba}), we finally obtain
%
\begin{eqnarray}
\sum_{k=0}^\infty\int_0^\infty \mathrm{e}^{-zt}
p_k^\nu(t) \,\mathrm{d}t & =& ( \lambda-\mu ) z^{\nu-1}
\int_0^\infty\frac{\mathrm{e}^{-sz^\nu}}{\lambda-
\mu \mathrm{e}^{-s ( \lambda-\mu )}}\nonumber\\
&&{} +
\int_0^\infty\frac{\mu- \mu \mathrm{e}^{-s ( \lambda- \mu )}}{
\lambda- \mu \mathrm{e}^{-s ( \lambda- \mu )}} z^{\nu-1} \mathrm{e}^{-s z^\nu}
\,\mathrm{d}s  \\
& = &\int_0^\infty z^{\nu-1} \mathrm{e}^{- s z^\nu} \,\mathrm{d}s = \frac{1}{z}, \nonumber
\end{eqnarray}
as desired.
\end{rem}

\begin{rem}
The distribution \eqref{pklgm} can be expressed in terms of the
probability law
of a fractional linear birth process with rate $\lambda-\mu$, which reads
\begin{eqnarray}
q_k^\nu ( t ) & =& \Pr \{ M_\nu ( t )
= k+l | M_\nu ( 0 ) = l+1 \}\nonumber\\[-8pt]\\[-8pt]
& =& \pmatrix{k+l-1\cr k-1} \sum_{r=0}^{k-1} ( -1 )^r
\pmatrix{k-1\cr r} E_{\nu,1} \bigl( - ( r+1+l ) ( \lambda
-\mu ) t^\nu \bigr), \nonumber
\end{eqnarray}
where $l+1$ initial progenitors are assumed (see \citet{pol}, formula (3.59)).
If we write
%
\begin{equation}
\Pr \{ G=l \} = \biggl( 1- \frac{\mu}{\lambda} \biggr)
\biggl( \frac{\mu}{\lambda} \biggr)^l,\qquad  l \geq0,
\end{equation}
then formula \eqref{pklgm} can be rewritten as\vspace*{1pt}
%
\begin{eqnarray}
\label{interpret}
p_k^\nu ( t ) & =& \biggl( \frac{\lambda-\mu}{\lambda} \biggr)^2
\sum_{l=0}^\infty\frac{l+k}{k} \biggl( \frac{\mu}{\lambda} \biggr)^l
\Pr \{ M_\nu(t) = k+l | M_\nu(0) =l+1 \} \nonumber\\
& =& \frac{\lambda-\mu}{\lambda}
\sum_{l=0}^\infty \biggl[ \biggl( 1+\frac{\mu}{k ( \lambda- \mu )}
\biggr) \Pr ( G=l ) +
\frac{\mu}{k} \frac{\mathrm{d}}{\mathrm{d} \mu} \Pr ( G=l ) \biggr]\\
&&\hphantom{\frac{\lambda-\mu}{\lambda}\sum_{l=0}^\infty}
{}\times\Pr \{ M_\nu ( t ) =
k+l | M_\nu ( 0 ) =
l+1 \} \nonumber
\end{eqnarray}
because
%
\begin{equation}
\frac{\mu}{k} \frac{\mathrm{d}}{\mathrm{d} \mu} \Pr(G=l) = \frac{l}{k}
\biggl( 1- \frac{\mu}{\lambda} \biggr) \biggl( \frac{\mu}{\lambda} \biggr)^l
- \frac{\mu}{k ( \lambda- \mu )} \Pr(G=l) .
\end{equation}
Result \eqref{interpret} shows that for large values of $k$, we have
the following,
interesting, approximation:
\[
p_k^\nu ( t ) {} \sim\frac{\lambda-\mu}{\lambda}
\sum_{l=0}^\infty\Pr ( G=l )
\Pr \{ M_\nu ( t ) =
k+l | M_\nu ( 0 ) =
l+1 \} .
\]
\end{rem}

\begin{thm}
For a fractional linear birth--death process $N_\nu ( t )$, $t>0$,
$\mu>\lambda$,
the probabilities $p_k^\nu ( t ) =
\Pr \{ N_\nu ( t ) = k \}$,
$k \geq1$, have
the following form:\vspace*{1pt}
\begin{eqnarray}
\label{probmgl}
p_k^\nu ( t ) &=& \biggl( \frac{\mu-\lambda}{\mu} \biggr)^2
\biggl( \frac{\lambda}{\mu} \biggr)^{k-1} \sum_{l=0}^\infty
\pmatrix{l+k\cr l} \biggl( \frac{\lambda}{\mu}
\biggr)^l\nonumber\\[-8pt]\\[-8pt]
&&\hphantom{\biggl( \frac{\mu-\lambda}{\mu} \biggr)^2
\biggl( \frac{\lambda}{\mu} \biggr)^{k-1} \sum_{l=0}^\infty}
{}\times\sum_{r=0}^{k-1} ( -1 )^r
\pmatrix{k-1\cr r}
E_{\nu,1} \bigl( - ( l+r+1 ) (\mu- \lambda
) t^\nu \bigr) .\nonumber
\end{eqnarray}
\end{thm}

\begin{pf}
By again using relationship \eqref{ite}, thanks to formula \eqref
{class-state} suitably rearranged, we can
write\vspace*{1pt}
%
\begin{equation}
p_k^\nu ( t ) = \int_0^\infty ( \mu-\lambda )^2
\mathrm{e}^{- ( \mu-\lambda )s} \lambda^{k-1}
\frac{ ( \mathrm{e}^{- ( \mu-\lambda )s} -1 )^{k-1}}{
( \lambda \mathrm{e}^{- ( \mu- \lambda )s} - \mu )^{k+1}}
\Pr \bigl( T_{2 \nu} ( t ) \in \mathrm{d}s \bigr) .
\end{equation}
By applying the Laplace transform, we have (omitting here some steps
similar to those of the proof of the previous theorem)\vspace*{1pt}
\begin{eqnarray}\label{lapmgl}
L_k^\nu ( z ) &= & \int_0^\infty
( \mu- \lambda )^2 \mathrm{e}^{- ( \mu-\lambda ) s}
( - \lambda )^{k-1} \frac{ ( 1-\mathrm{e}^{- (
\mu-\lambda ) s} )^{k-1}}{ ( -\mu )^{k+1}
( 1 - (\lambda/\mu) \mathrm{e}^{- ( \mu- \lambda ) s}
)^{k+1}} z^{\nu-1} \mathrm{e}^{-s z^\nu} \,\mathrm{d}s \nonumber\\
&=  & \biggl( \frac{\mu-\lambda}{\mu} \biggr)^2 \biggl( \frac{\lambda}{
\mu} \biggr)^{k-1} \sum_{l=0}^\infty
\pmatrix{l+k\cr l}\biggl ( \frac{\lambda}{\mu}
\biggr)^l
\sum_{r=0}^{k-1} ( -1 )^r \pmatrix{k-1\cr r}
z^{\nu-1}\\ 
&&\hphantom{\biggl( \frac{\mu-\lambda}{\mu} \biggr)^2 \biggl( \frac{\lambda}{
\mu} \biggr)^{k-1} \sum_{l=0}^\infty
\pmatrix{l+k\cr l}\biggl ( \frac{\lambda}{\mu}
\biggr)^l
\sum_{r=0}^{k-1}}
{}\times \int_0^\infty \mathrm{e}^{-s ( z^\nu+ ( \mu- \lambda )
( l+r+1 ) ) } \,\mathrm{d}s \nonumber\\
&=  & \biggl( \frac{\mu-\lambda}{\mu} \biggr)^2 \biggl( \frac{\lambda}{
\mu} \biggr)^{k-1} \sum_{l=0}^\infty
\pmatrix{l+k\cr l} \biggl( \frac{\lambda}{\mu}
\biggr)^l
\sum_{r=0}^{k-1} ( -1 )^r \pmatrix{k-1\cr r}
\frac{z^{\nu-1}}{z^\nu+ ( \mu-\lambda )
( l+r+1 )} \nonumber.
\end{eqnarray}
By transforming equation \eqref{lapmgl}, we easily arrive at the result
\eqref{probmgl}.
\end{pf}

\begin{rem}
When $k=1$, equation \eqref{pklgm} takes a simple form:
\begin{eqnarray}
p_1^\nu ( t ) & =&
\biggl( \frac{\lambda-\mu}{\lambda} \biggr)^2 \sum_{l=0}^\infty
( l+1 ) \biggl( \frac{\mu}{\lambda} \biggr)^l
E_{\nu,1} \bigl( - ( l+1 ) ( \lambda-\mu ) t^\nu \bigr)
\nonumber\\[-8pt]\\[-8pt]
& = &\biggl( \frac{\lambda-\mu}{\lambda} \biggr)^2 \sum_{l=1}^\infty
l \biggl( \frac{\mu}{\lambda} \biggr)^{l-1}
E_{\nu,1} \bigl( - l ( \lambda-\mu ) t^\nu \bigr)
\nonumber,
\end{eqnarray}
where $\lambda>\mu$. For the case $\lambda<\mu$, we obtain
essentially the same expression with $\lambda$ and $\mu$
exchanged.
\end{rem}

An interpretation similar to that in \eqref{interpret} is valid for the
case $\mu>\lambda$ as well.
The following theorem describes the structure of the state probabilities
$p_k^\nu ( t )$, $k \geq1$, in the case where $\mu=\lambda$, that is,
when the birth rate equals the
death rate.
\begin{thm}
In the case $\mu=\lambda$, the probabilities $p_k^\nu ( t ) =
\Pr \{ N_\nu ( t ) = k \}$ of the fractional
linear birth--death process read
%
\begin{equation}
\label{formulaC}
\Pr \{ N_\nu ( t ) = k \} = \frac{
( -1 )^{k-1}
\lambda^{k-1}}{k!} \frac{\mathrm{d}^k}{\mathrm{d} \lambda^k} \bigl[ \lambda \bigl( 1-
p_0^\nu ( t ) \bigr) \bigr]
\end{equation}
with $k \geq1$ and $t>0$.
\end{thm}

\begin{pf}
The explicit form of the distribution $\Pr \{ N_\nu (
t ) = k \},
k \geq1$, for the fractional linear birth--death process, in
the case $\lambda=\mu$, can be evaluated in
the following manner.
In light of \eqref{sub-repr}, we have
%
\begin{equation}
\Pr \{ N_\nu ( t ) = k \}
= \int_0^\infty\Pr
\{ N ( s ) = k \} \Pr \{ T_{2 \nu} ( t
) \in \mathrm{d}s \}
\end{equation}
so that
%
\begin{eqnarray}
\label{rel-dist}
L_k^\nu ( z ) &=& \int_0^\infty \mathrm{e}^{- z t} \Pr \{
N_\nu ( t ) =k \} \,\mathrm{d}t\nonumber\\[-8pt]\\[-8pt]
& = &\int_0^\infty
\frac{ ( \lambda s
)^{k-1}}{ ( 1+ \lambda s )^{k+1} }
z^{\nu-1} \mathrm{e}^{-s z^\nu} \,\mathrm{d}s .\nonumber
\end{eqnarray}
This is because for the $\lambda=\mu$ case of the classical
birth--death process,
we have that
(see \citet{bailey}, formula (8.53), page 95)
%
\begin{equation}
\label{formulaD}
\Pr \{ N ( t ) = k \} =
\frac{ ( \lambda t )^{
k-1}}{ ( 1+\lambda t )^{k+1}},\qquad  k \geq1 .
\end{equation}
We note that the extinction probability cannot be extracted
from the above formula since
it reads
%
\begin{equation}
\Pr \{ N ( t ) =0 \} =
\frac{\lambda t}{1+ \lambda t} .
\end{equation}
This implies that we have a different expression for $k \geq1$ and
$k=0$ for the fractional linear
birth--death process as well.

Formula \eqref{rel-dist} can be expanded out as
\begin{eqnarray}
\label{formulaA}
L_k^\nu ( z ) & =& \frac{ ( -1 )^k \lambda^{k-1}}{k!}
\frac{\mathrm{d}^k}{\mathrm{d} \lambda^k} \int_0^\infty\frac{1}{s (
1 + \lambda s )}
z^{\nu-1} \mathrm{e}^{-s z^\nu} \,\mathrm{d}s \nonumber\\
& =& \frac{ ( -1 )^k \lambda^{k-1}}{k!} \frac{\mathrm{d}^k}{\mathrm{d} \lambda^k}
\int_0^\infty \biggl( \frac{1}{s} - \frac{\lambda}{1+ \lambda
s} \biggr) z^{\nu-1}
\mathrm{e}^{-s z^\nu} \,\mathrm{d}s\nonumber\\[-8pt]\\[-8pt]
& = &\frac{ ( -1 )^k \lambda^{k-1}}{k!} \frac{\mathrm{d}^k}{\mathrm{d} \lambda^k}
\biggl[ \int_0^\infty\int_0^\infty ( \mathrm{e}^{-w s} - \lambda \mathrm{e}^{-w (
1+ \lambda s )}
) z^{\nu-1} \mathrm{e}^{-s z^\nu} \,\mathrm{d}s\, \mathrm{d}w \biggr] \nonumber\\
& =& \frac{ ( -1 )^k \lambda^{k-1}}{k!} z^{\nu-1}
\frac{\mathrm{d}^k}{\mathrm{d} \lambda^k}
\biggl[ - \int_0^\infty\frac{ \lambda \mathrm{e}^{-w}}{w \lambda+ z^\nu}
\,\mathrm{d}w + \int_0^\infty
\frac{\mathrm{d}w}{w + z^\nu} \biggr] . \nonumber
\end{eqnarray}
By inverting the Laplace transform, we have that
%
\begin{eqnarray}
\label{formulaB}
\Pr \{ N_\nu ( t ) = k \} & =&
\frac{ ( -1 )^k \lambda^{k-1}}{k!} \frac{\mathrm{d}^k}{\mathrm{d} \lambda^k}
\biggl[ \int_0^\infty \bigl( E_{\nu,1} ( - w t^\nu )
- \lambda \mathrm{e}^{-w} E_{\nu,1} ( -
\lambda w t^\nu
) \bigr)\, \mathrm{d}w \biggr]\nonumber
\\
& = &\frac{ ( -1 )^{k-1} \lambda^{k-1}}{k!}
\frac{\mathrm{d}^k}{\mathrm{d} \lambda^k}
\biggl[ \lambda\int_0^\infty \mathrm{e}^{-w} E_{\nu,1} (
- \lambda w t^\nu )
\,\mathrm{d}w \biggr] \\
& = &\frac{ ( -1 )^{k-1} \lambda^{k-1}}{k!}
\frac{\mathrm{d}^k}{\mathrm{d} \lambda^k}
\bigl[ \lambda \bigl( 1 - p_0^\nu ( t ) \bigr)
\bigr] . \nonumber
\end{eqnarray}
Formula \eqref{formulaC} is thus proved.
\end{pf}

It is important to note how all the state probabilities $p_k^\nu
( t )$ depend
on the extinction probability $p_0^\nu ( t )$.

\begin{rem}
For $\nu=1$, we can extract from \eqref{formulaC} the classical formula
\eqref{formulaD} because
%
\begin{equation}
p_k^1 ( t ) = \Pr \{ N ( t ) = k \} =
\frac{ ( -1 )^{k-1} \lambda^{k-1}}{k!} \frac{\mathrm{d}^k}{\mathrm{d} \lambda^k}
\biggl[ \frac{\lambda}{1+ \lambda t} \biggr]
\end{equation}
and because
%
\begin{eqnarray}
\label{considering}
\frac{\mathrm{d}^k}{\mathrm{d} \lambda^k} \biggl[ \frac{\lambda}{1+\lambda t} \biggr] & =&
\sum_{j=0}^k \pmatrix{k\cr j} \frac{\mathrm{d}^j}{\mathrm{d} \lambda^j} \lambda
\frac{\mathrm{d}^{k-j}}{\mathrm{d} \lambda^{k-j}} \biggl( \frac{1}{1+\lambda t} \biggr) \nonumber\\
& = &\lambda\frac{\mathrm{d}^k}{\mathrm{d} \lambda^k} \biggl( \frac{1}{1+\lambda t} \biggr)
+ k \frac{\mathrm{d}^{k-1}}{\mathrm{d} \lambda^{k-1}} \biggl( \frac{1}{1+\lambda t} \biggr)
\nonumber\\
& =& \lambda\frac{ ( -1 )^k k! t^k}{ ( 1+ \lambda t )^{k+1}}
+ k \frac{ ( -1 )^{k-1} ( k-1 )! t^{k-1}}{
( 1+\lambda t )^k} \\
& =& \frac{ ( k-1 )! t^{k-1}}{ ( 1+\lambda t )^{k+1}}
( -1 )^{k-1} [ - \lambda k t+ ( 1+\lambda t ) k
] \nonumber\\
& =& \frac{k! t^{k-1}}{ ( 1+\lambda t )^{k+1}} ( -1 )^{k-1}
. \nonumber
\end{eqnarray}
\end{rem}

\begin{rem}
From the representation on the last line of \eqref{formulaB}, it is
possible to give an alternative proof of the subordination relationship
\eqref{ite} when $k \geq1$, as follows:
%
\begin{eqnarray}
p_k^\nu ( t ) & = &\frac{ ( -1 )^{k-1} \lambda^{k-1}}{
k!} \frac{\mathrm{d}^k}{\mathrm{d} \lambda^k} \bigl[ \lambda \bigl( 1 - p_0^\nu ( t )
\bigr) \bigr] \\
& =& \frac{ ( -1 )^{k-1} \lambda^{k-1}}{
k!} \frac{\mathrm{d}^k}{\mathrm{d} \lambda^k} \biggl[ \lambda- \int_0^\infty
\frac{\lambda^2 s}{1+ \lambda s} \Pr \bigl( T_{2 \nu} ( t )
\in \mathrm{d}s \bigr) \biggr] \nonumber\\
& =& \frac{ ( -1 )^{k-1} \lambda^{k-1}}{
k!} \biggl[ \int_0^\infty\frac{\mathrm{d}^k}{\mathrm{d} \lambda^k} \biggl[ \frac{\lambda}{
1+\lambda s} \biggr] \Pr \bigl( T_{2 \nu} ( t ) \in \mathrm{d}s
\bigr) \biggr] .
\end{eqnarray}
Exploiting \eqref{considering}, we readily obtain
%
\begin{equation}
p_k^\nu ( t ) = \int_0^\infty\frac{ ( \lambda s )^{k-1}}{
( 1+\lambda s )^{k+1} } \Pr \bigl( T_{2 \nu} ( t
) \in \mathrm{d}s \bigr) .
\end{equation}
\end{rem}

\begin{rem}
Here, we present two other interesting relationships. The first one is simply
a particular case of formula \eqref{formulaC} when $k=1$, that is, the
probability of having
one individual in the process at time $t$ is
%
\begin{equation}
\Pr \{ N_\nu ( t ) = 1 \} =
\frac{\mathrm{d}}{\mathrm{d} \lambda} \bigl[ \lambda \bigl( 1- p_0^\nu ( t )
\bigr) \bigr] .
\end{equation}
The second relationship is again a particular case of formula
\eqref{formulaC} with $\nu=1/2$. In that case, recalling that
%
\begin{equation}
E_{{1}/{2},1} ( x ) = \frac{2}{\sqrt{\uppi}}
\int_0^\infty \mathrm{e}^{- y^2 + 2 y x}\, \mathrm{d}y,
\end{equation}
we obtain
%
\begin{eqnarray}
\Pr \{ N_{{1}/{2}} ( t ) = k \} & =&
\frac{ ( -1 )^{k-1} \lambda^{k-1}}{k!}
\frac{\mathrm{d}^k}{\mathrm{d} \lambda^k} \bigl[ \lambda \bigl( 1- p_0^{1/{2}} ( t
) \bigr) \bigr] \nonumber\\
& = &\frac{ ( -1 )^{k-1} \lambda^{k-1}}{k!}
\frac{\mathrm{d}^k}{\mathrm{d} \lambda^k} \biggl[ \lambda\int_0^\infty \mathrm{e}^{-w} E_{1/2,1}
( - \lambda t^{1/{2}} w ) \,\mathrm{d}w \biggr] \nonumber\\
& =& \frac{ ( -1 )^{k-1} \lambda^{k-1}}{k!}\nonumber\\
&&{}\times\frac{\mathrm{d}^k}{\mathrm{d} \lambda^k} \biggl[ \frac{2 \lambda}{\sqrt{\uppi}} \int_0^\infty
\mathrm{e}^{-w} \int_0^\infty \mathrm{e}^{-y^2 - 2 y \lambda t^{{1}/{2}} w} \,\mathrm{d}w \,\mathrm{d}y \biggr]\qquad
\\
& =& \frac{ ( -1 )^{k-1} \lambda^{k-1}}{k!}
\frac{\mathrm{d}^k}{\mathrm{d} \lambda^k} \biggl[ \frac{2 \lambda}{\sqrt{\uppi}} \int_0^\infty
\frac{\mathrm{e}^{-y^2}}{1+2 \lambda y \sqrt{t}} \,\mathrm{d}y \biggr] \nonumber\\
& =& \frac{ ( -1 )^{k-1} \lambda^{k-1}}{k!}
\frac{\mathrm{d}^k}{\mathrm{d} \lambda^k} \biggl[ 2 \lambda\int_0^\infty
\frac{\mathrm{e}^{-w^2/{(2 t)}}}{1 + \lambda\sqrt{2} w} \frac{1}{
\sqrt{2 \uppi t}} \,\mathrm{d}w \biggr] \nonumber\\
& = &\frac{ ( -1 )^{k-1} \lambda^{k-1}}{k!}
\frac{\mathrm{d}^k}{\mathrm{d} \lambda^k} \mathbb{E} \biggl[ \frac{ 2 \lambda}{
1 + \lambda\sqrt{2} B ( t )} \biggr], \nonumber
\end{eqnarray}
where $B ( t )$, $t>0$, is a standard Brownian motion.
\end{rem}

\section{Some further properties of the fractional linear birth--death process}
The analysis of the moments of the fractional linear birth--death
process gives us useful information
concerning the behaviour of the system.
Starting from \eqref{eq:fractDiffEqGF}, we easily see that
%
\begin{equation}
\label{eq:mean-value}
\mathbb{E} N_\nu ( t ) =
\frac{\partial G}{\partial u} \bigg|_{u=1}
\end{equation}
is the solution to
%
\begin{equation}
\label{eq:eq-mean-value-solution}
\cases{
\dfrac{ \mathrm{d}^\nu}{\mathrm{d} t^\nu} \mathbb{E} N_\nu= ( \lambda- \mu )
\mathbb{E} N_\nu, & \quad $0<\nu\leq1$,\vspace*{2pt} \cr
\mathbb{E} N_\nu ( 0 ) = 1.
}
\end{equation}
By again applying the Laplace transform, we have that the solution
to \eqref{eq:eq-mean-value-solution}
is
%
\begin{equation}
\label{eq:result-expect}
\mathbb{E} N_\nu ( t ) = E_{\nu,1} \bigl( (
\lambda- \mu ) t^\nu \bigr),\qquad  t>0 .
\end{equation}
In the case $\lambda> \mu$, the result \eqref{eq:result-expect} shows
that the mean size
of the population coincides with that of a fractional linear pure
birth process
with rate $\lambda- \mu> 0$ (see \citet{pol}).
Result \eqref{eq:result-expect} can also be derived by means of
the subordination relationship
\eqref{ite}:
%
\begin{eqnarray}
\label{eq:other-expect}
\mathbb{E} N_\nu ( t ) & =& \sum_{k=0}^\infty k
\Pr \{ N_\nu (
t ) = k \} \nonumber\\
& =& \sum_{k=0}^\infty k \int_0^\infty\Pr \{ N
( s ) = k \}
\Pr \{ T_{2 \nu} ( t ) \in \mathrm{d}s \} \\
& =& \int_0^\infty \mathrm{e}^{ ( \lambda- \mu )
s} \Pr \{ T_{2 \nu} ( t )
\in \mathrm{d}s \} . \nonumber
\end{eqnarray}
The Laplace transform of (\ref{eq:other-expect}) yields
%
\begin{eqnarray}
\int_0^\infty \mathrm{e}^{- z t} \mathbb{E} N_\nu (
t ) \,\mathrm{d}t & = &\int_0^\infty \mathrm{e}^{ (\lambda- \mu
) s} z^{\nu-1} \mathrm{e}^{- s z^\nu} \,\mathrm{d}s \\
& =& \frac{z^{\nu-1}}{z^\nu- ( \lambda-
\mu ) } = \int_0^\infty \mathrm{e}^{ - z t} E_{\nu,1}
\bigl( ( \lambda- \mu )
t^\nu \bigr)\, \mathrm{d}t \nonumber
\end{eqnarray}
and this confirms \eqref{eq:result-expect}.

By again applying \eqref{eq:fractDiffEqGF}, it is also possible to derive
the variance $\mathbb{V} $ar$N_\nu ( t )$, $t>0$, of the number of
individuals in the population at time $t$. We start by evaluating the
second-order
factorial moment $\mu_{ ( 2 )} ( t ) = \mathbb{E}
[ N_\nu ( t ) ( N_\nu ( t ) -1 ) ]$, $t>0$.
From \eqref{eq:fractDiffEqGF}, after some straightforward steps, we
see that
%
\begin{equation}
\mu_{ ( 2 )} ( t ) = \mathbb{E}
\bigl[ N_\nu ( t ) \bigl( N_\nu ( t ) -1 \bigr) \bigr]
= \frac{\partial^2 G}{\partial u^2} \bigg|_{u=1}
\end{equation}
is the solution to the following differential equation:
%
\begin{equation}
\label{momsec}
\cases{
\dfrac{\mathrm{d}^\nu}{\mathrm{d} t^\nu} \mu_{ ( 2 )} ( t )
= 2 \lambda\mathbb{E} N_\nu ( t ) + 2 ( \lambda-\mu )
\mu_{ ( 2 )} ( t ) , &\quad $ 0<\nu\leq1$, \cr
\mu_{ ( 2 )} ( 0 ) = 0.
}
\end{equation}
In order to solve \eqref{momsec}, we apply the Laplace transform,
obtaining, in the case
$\lambda\neq\mu$,
\begin{eqnarray}
\label{lapsec}
\int_0^\infty \mathrm{e}^{-zt} \mu_{ ( 2 )} ( t )\,\mathrm{d}t & =&
2 \lambda\frac{z^{\nu-1}}{z^\nu- ( \lambda- \mu )} \cdot
\frac{1}{z^\nu- 2 ( \lambda- \mu )} \nonumber\\[-8pt]\\[-8pt]
& =&\frac{2 \lambda z^{\nu-1}}{\lambda- \mu} \biggl[ \frac{1}{z^\nu-2 (
\lambda- \mu )} - \frac{1}{z^\nu- ( \lambda- \mu )} \biggr] .
\nonumber
\end{eqnarray}
The Laplace transform \eqref{lapsec} can be inverted, thus leading to
the explicit
expression of the second-order factorial moment as
%
\begin{equation}
\label{secord}
\mu_{ ( 2 )} ( t ) = \frac{2 \lambda}{\lambda-\mu}
\bigl[ E_{\nu,1} \bigl( 2 ( \lambda-\mu ) t^\nu \bigr) -
E_{\nu,1} \bigl( ( \lambda-\mu ) t^\nu \bigr) \bigr] .
\end{equation}
From the first expression of the Laplace transform in \eqref{lapsec},
we also have that
%
\begin{equation}
\mu_{ ( 2 )} ( t ) = 2 \lambda\int_0^t s^{\nu-1}
E_{\nu,\nu} \bigl( 2 ( \lambda- \mu ) s^\nu \bigr) E_{\nu,1}
\bigl( ( \lambda- \mu ) ( t-s )^\nu \bigr)\, \mathrm{d}s .
\end{equation}
By applying similar calculations to those of \eqref{conv}, we prove result
\eqref{secord}.

From \eqref{secord}, we can easily write that
%
\begin{eqnarray}
\label{resultsec}
\operatorname{\mathbb{V}ar}N_\nu ( t ) & = & \frac{2 \lambda}{\lambda- \mu}
\bigl[ E_{\nu,1} \bigl( 2 ( \lambda- \mu ) t^\nu \bigr) -
E_{\nu,1} \bigl( ( \lambda-\mu ) t^\nu \bigr) \bigr] \nonumber \\
&&{} + E_{\nu,1} \bigl( ( \lambda-\mu ) t^\nu \bigr) -
E_{\nu,1}^2 \bigl( ( \lambda-\mu ) t^\nu \bigr) \\
& = & \frac{2 \lambda}{\lambda-\mu} E_{\nu,1} \bigl( 2 ( \lambda- \mu
) t^\nu \bigr) - \frac{\lambda+ \mu}{\lambda-\mu} E_{\nu,1}
\bigl( ( \lambda- \mu ) t^\nu \bigr) - E_{\nu,1}^2 \bigl(
( \lambda-\mu ) t^\nu \bigr) \nonumber.
\end{eqnarray}
\begin{rem}
When $\nu=1$, we obtain from \eqref{resultsec} the expression for the variance
of the classical linear birth--death process as follows:
\begin{eqnarray}
\label{clasvar}
\operatorname{\mathbb{V}ar}N ( t ) & =& \frac{2 \lambda}{
\lambda-\mu} \mathrm{e}^{2 t ( \lambda- \mu )} -
\frac{\lambda+ \mu}{\lambda- \mu} \mathrm{e}^{t ( \lambda- \mu )}
- \mathrm{e}^{2 t ( \lambda- \mu )} \nonumber\\[-8pt]\\[-8pt]
& =& \frac{\lambda+ \mu}{\lambda- \mu} \bigl( \mathrm{e}^{2 t (
\lambda- \mu )} - \mathrm{e}^{t ( \lambda- \mu )} \bigr) =
\frac{\lambda+\mu}{\lambda- \mu} \mathrm{e}^{t ( \lambda- \mu )}
\bigl( \mathrm{e}^{t ( \lambda- \mu )} -1 \bigr) . \nonumber
\end{eqnarray}
\end{rem}
\begin{rem}
When $\mu=0$, that is, in the case of pure linear birth, we obtain
from \eqref{clasvar}
the expression of the variance of the classical linear
pure birth process and from \eqref{resultsec}
that of the fractional linear birth process (see \citet{pol}).
\end{rem}

In the case $\lambda=\mu$, from \eqref{momsec}, it is easy to show that
%
\begin{equation}
\operatorname{\mathbb{V}ar}N_\nu ( t ) = \frac{2 \lambda t^\nu}{\Gamma (
\nu+1 )},
\end{equation}
in accordance with the well-known result of the classical linear
birth--death process
for $\lambda=\mu$ which reads $\operatorname{\mathbb{V}ar}N ( t ) = 2 \lambda t$.
\begin{rem}
We can directly evaluate the mean value $\mathbb{E} N ( t )$ for
$\lambda=\mu$ in the following way:
%
\begin{eqnarray}
\mathbb{E} N ( t ) & =& \sum_{k=1}^\infty k \biggl( \frac{
( \lambda t )^{k-1}}{ ( 1+ \lambda t )^{k+1}}
\biggr)
= \frac{1}{ ( 1 + \lambda t )^2} \sum_{k=1}^\infty
k \biggl( \frac{\lambda t}{1+ \lambda t} \biggr)^{k-1} \nonumber\\
& =& \frac{1}{ ( 1 + \lambda t )^2} \frac{\mathrm{d}}{\mathrm{d}z}
\sum_{k=1}^\infty z^k \Bigg|_{z={\lambda t}/{(1+ \lambda t)}}
= \frac{1}{ ( 1 + \lambda t )^2} \frac{\mathrm{d}}{\mathrm{d}z}
\frac{z}{1-z} \bigg|_{z={\lambda t}/{(1 + \lambda t)}} \\
& =& \frac{1}{ ( 1+ \lambda t )^2} \frac{1}{ ( 1-z )^2}
\bigg|_{z={\lambda t}/{(1 + \lambda t)}} =1 . \nonumber
\end{eqnarray}
The assumption that $\lambda=\mu$ implies that the mean size of the population
$\mathbb{E} N_\nu(t)$, $t>0$, is equal to one (number of original progenitors)
for all $t>0$ and all
$0<\nu\leq1$ (this is also confirmed
for $\lambda=\mu$ by \eqref{eq:result-expect}).
\end{rem}

\section*{Acknowledgements}
The authors wish to thank Francis Farrelly
for having
checked and corrected the manuscript. Thanks are due to the referee
who detected
misprints and flaws.

\printhistory

\end{document}